\documentclass[%aip,
jmp,
amsmath,amssymb,
reprint,%
numerical,
]{revtex4-1}
\usepackage{graphicx}% Include figure files
\usepackage{dcolumn}% Align table columns on decimal point
\usepackage{bm}% bold math
\usepackage{hyperref}
\usepackage[utf8]{inputenc}
\usepackage[T1]{fontenc}
\usepackage{mathptmx}
\usepackage{etoolbox}
\newcommand{\RNum}[1]{\uppercase\expandafter{\romannumeral #1\relax}}
\usepackage{array,xcolor,colortbl}
\usepackage{multirow}
\newcommand{\ud}{\ensuremath{\mathrm{d}}}

\newcommand{\beq}{\begin{equation}}
\newcommand{\eeq}{\end{equation}}

\makeatletter
\def\@email#1#2{%
	\endgroup
	\patchcmd{\titleblock@produce}
	{\frontmatter@RRAPformat}
	{\frontmatter@RRAPformat{\produce@RRAP{*#1\href{mailto:#2}{#2}}}\frontmatter@RRAPformat}
	{}{}
}%
\makeatother

\begin{document}
	
\title{Thermal phototactic bioconvection in an isotropic porous medium heated from above}
% Force line breaks with \\

% Force line breaks with \\
\author{S. K. Rajput}
\altaffiliation[Corresponding author: E-mail: ]{shubh.iiitj@gmail.com}%Lines break automatically or can be forced with \\
\author{M. K. Panda}%
%\email{mkpanda@iiitdmj.ac.in}
\author{A. Rathi}%
%\email{shubh.iiitj@gmail.com}
\affiliation{Department of Mathematics, PDPM Indian Institute of Information Technology Design and Manufacturing, Jabalpur 482005, India.%\\This line break forced with \textbackslash\textbackslash
}%
%\homepage{http://www.Second.institution.edu/~Charlie.Author.}
%\affiliation{%
%	Second institution and/or address%\\This line break forced% with \\
%}%

\begin{abstract}	
	This study investigates thermal phototactic bioconvection in an isotropic porous medium using the Darcy-Brinkman model. The top boundary of the medium is exposed to normal collimated light and subjected to heating. A linear analysis of bio-thermal convection is performed using a fourth-order accurate finite difference scheme, employing Newton-Raphson-Kantorovich iterations for both rigid-free and rigid-rigid boundary conditions. The effects of the Lewis number, Darcy number, and thermal Rayleigh number on bioconvective processes are examined and presented graphically. The findings reveal that increasing the thermal Rayleigh number stabilizes the suspension, whereas a higher Lewis number enhances instability.
\end{abstract}

%%%%%%%%%%%%%%%%%%%%%%%%%%%%%%%%%%%%%%%%%%%%%%%%%%%%%%%%%%%%%%%%%%%%%%%%%%%	

\maketitle

%%%%%%%%%%%%%%%%%%%%%%%%%%%%%%%%%%%%%%%%%%%%%%%%%%%%%%%%%%%%%%%%%%%%%%%%%%%	
\section{INTRODUCTION}

Bioconvection is a fluid dynamics phenomenon driven by the collective motion of self-propelled microorganisms that are slightly denser than their surrounding fluid. This process plays a crucial role in various biological and industrial applications~\cite{19platt1961,20pedley1992,21hill2005,22bees2020,23javadi2020}. Motile microorganisms, such as algae and bacteria, generate density variations as they respond to external stimuli, a behavior known as taxis, leading to convective instability. Key examples of taxis include phototaxis, gravitaxis, gyrotaxis, chemotaxis, and thermotaxis. Understanding bioconvection is particularly relevant in environmental science, biotechnology, and engineering, where it influences nutrient transport, bioreactor efficiency, and microbial ecology.

Early research primarily focused on suspensions under isothermal conditions. However, many microorganisms, particularly thermophiles inhabiting hot springs, thrive in environments with significant temperature variations~\cite{27kuznetsov2005thermo,28alloui2006stability,29nield2006onset}. Among the various types of taxis influencing microbial movement, phototaxis (response to light) and thermotaxis (response to temperature gradients) play a crucial role in shaping bioconvective patterns~\cite{39rajput2025d}. While substantial research has been dedicated to phototactic and gravitactic bioconvection in non-porous media, thermal bioconvection in porous media saturated with algal suspensions remains relatively under-explored. The presence of a porous matrix introduces additional complexities, such as flow resistance and modified convective dynamics, making it a critical area of study for both natural ecosystems and industrial fluid systems.

\begin{figure}[!htbp]
	\centering
	\includegraphics[scale=0.5]{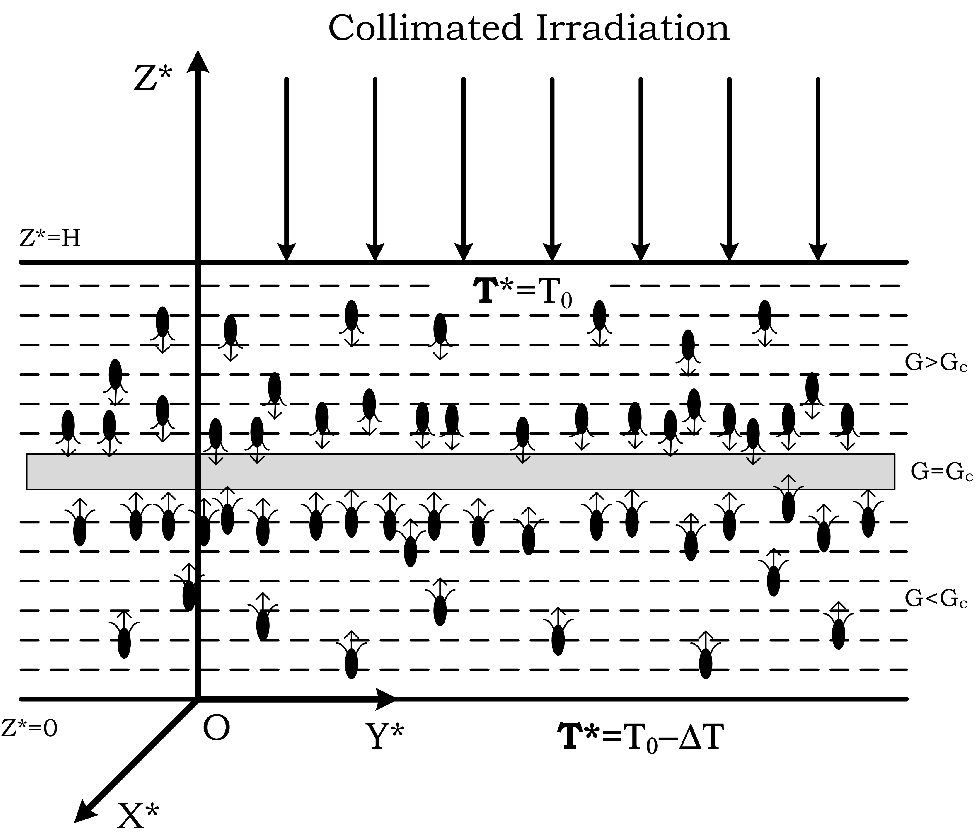}
	\caption{\footnotesize{The spatial configuration of the proposed problem.}}
	\label{fig1}
\end{figure}
	
The formation and characteristics of bioconvection patterns influenced by phototaxis depend on various environmental light conditions, including direct and oblique collimated solar irradiation~\cite{1wager1911,2kitsunezaki2007,3kessler1985,4williams2011,24kessler1986,25kessler1997,12vincent1996}. High-intensity light can either disrupt established patterns or inhibit their development~\cite{3kessler1985,4williams2011,5kessler1989}. Variations in illumination levels contribute to changes in the spatial structure and size of these patterns. These alterations can be attributed to specific mechanisms. Photosynthetic microorganisms exhibit directional movement in response to light intensity. When the intensity $G$ remains below a critical threshold $G_c$, they exhibit positive phototaxis, migrating toward brighter regions. However, when $G$ surpasses $G_c$, they exhibit negative phototaxis, moving toward lower light intensities to avoid potential damage. This behavior results in the accumulation of microorganisms in regions where $G = G_c$, forming structured bioconvective patterns (see Fig.~\ref{fig1})~\cite{6hader1987}. Additionally, the process of light absorption significantly influences pattern formation by affecting the spatial distribution of microorganisms within the suspension~\cite{7ghorai2010,9straughan1993}.  

Over the years, significant progress has been made in understanding the mechanisms of bioconvection. Early research primarily focused on gravitactic microorganisms, where the swimming behavior of dense motile microbes was analyzed to determine the instability conditions that lead to convective flow~\cite{5kessler1989}. Later, studies on phototactic bioconvection emerged, emphasizing the interactions between microbial responses to light intensity, absorption effects, and microorganism concentration~\cite{12vincent1996}. Ghorai and Hill~\cite{10ghorai2005} extended this research by developing two-dimensional phototactic bioconvection models, emphasizing the role of shading effects. Further studies by Ghorai \textit{et al.}~\cite{7ghorai2010} investigated phototactic bioconvection in isotropically scattering algal suspensions, while Panda and Singh~\cite{11panda2016} analyzed its onset in a two-dimensional system with free side walls. The scattering properties of algae, influenced by their shape and size, affect light propagation. Early studies suggested that algae predominantly scatter light in the forward direction. To explore the role of forward scattering at the onset of bioconvection, Ghorai and Panda~\cite{13ghorai2013} developed a mathematical model. Panda \textit{et al.}~\cite{15panda2016} later investigated the combined effects of diffuse and collimated light fluxes in an isotropic scattering medium, while Panda~\cite{8panda2020} extended this research to anisotropic scattering media. Further work by Panda \textit{et al.}~\cite{16panda2022} examined the impact of oblique collimated flux on bioconvective instability in algal suspensions, whereas Panda and Rajput~\cite{26rajput2023} studied the combined effects of diffuse and oblique collimated flux in a uniformly (isotropic) scattering suspension. Recent studies by Rajput and Panda~\cite{35rajput2024,37rajput2025b} focused on the influence of a diffuse flux alone on photo-bioconvective instability in an isotropically scattering suspension with both rigid and free top walls. More recently, the same authors~\cite{36rajput2025a,38rajput2025c} analyzed bioconvective instability in an anisotropic suspension within the same geometric framework.

In addition to phototaxis, thermal effects on bioconvection have attracted considerable interest due to their influence on microorganism motility, metabolic activity, and fluid density variations. Kuznetsov~\cite{40kuznetsov2006} investigated thermo-bioconvection in a shallow porous layer heated from below, highlighting the role of thermal Rayleigh numbers in convective instability. Sheremet and Pop~\cite{41sheremet2014thermo} examined thermo-bioconvection within a square porous cavity, focusing on the behavior of oxytactic microorganisms in the presence of thermal gradients. Zhao et al.~\cite{42zhao2019} explored the impact of porous matrices on thermal bioconvection, providing insights into heat transfer mechanisms governing bioconvective flows. Further studies have examined external influences on thermal bioconvection. Biswas et al.~\cite{44biswas2021thermo} analyzed the role of magnetic fields in thermo-bioconvection, proposing new strategies for regulating convection in bio-suspensions. More recently, Kopp and Yanovsky~\cite{45kopp2024weakly} incorporated weakly nonlinear stability analysis to investigate bio-thermal convection under rotation, gravity modulation, and heat sources, contributing to a deeper theoretical understanding of the phenomenon.

Despite significant advancements, the intricate interplay between thermal effects, porous media, and phototaxis in bioconvection remains under-explored. This study seeks to bridge these gaps by formulating a mathematical model for thermal bioconvection in a porous medium. The model incorporates the Navier-Stokes equations with the Boussinesq approximation, a microorganism conservation equation, and an energy equation. The primary objectives include analysing stability conditions, examining the influence of various physical parameters on convection onset, and assessing practical implications for biotechnological applications. These applications span diverse fields, including biofuel production, water filtration, and environmental remediation.

This paper is structured as follows: Section 2 introduces the mathematical formulation, defining the governing equations and boundary conditions. Section 3 discusses the steady-state solution, outlining the equilibrium conditions. Section 4 presents the linear stability analysis, which establishes the criteria for instability. Section 5 describes the numerical solution methodology, followed by Section 6, which provides numerical results and examines the influence of various parameters on bioconvection behavior. Finally, Section 7 summarizes the key findings and proposes potential directions for future research.

%%%%%%%%%%%%%%%%%%%%%%%%%%%%%%%%%%%%%%%%%%%%%%%%%%%%%%%%%%%%%%%%%%%%	

\section{Problem formulation}
A porous medium filled with an algal suspension is considered, confined between two horizontal planes located at \( z^* = 0 \) and \( z^* = H \), extending infinitely in the horizontal directions. The study examines two boundary configurations: (i) a free-rigid setup, where the upper boundary is stress-free and the lower boundary is rigid, and (ii) a rigid-rigid setup, where both boundaries are rigid. The system is illuminated by collimated irradiation from above, while thermal effects are incorporated by applying heat at the bottom boundary. To ensure microbial viability and preserve phototactic behavior, the imposed temperature variations are assumed to be minimal. The upper boundary is maintained at a constant temperature \( T_0 \), while the lower boundary is kept at \( T_0 - \Delta T \).

%%%%%%%%%%%%%%%%%%%%%%%%%%%%%%%%%%%%%%%%%%%%%%%%%%%%%%%%%%%%%%%%%%	     
\subsection{THE GOVERNING EQUATIONS}
A dilute, incompressible suspension of phototactic microorganisms is considered. The governing equations are formulated based on the models introduced by Kuznetsov~\cite{27kuznetsov2005thermo} and Vincent and Hill~\cite{12vincent1996}. It is assumed that the heating process does not alter the microorganisms' phototactic behavior, including their swimming orientation and velocity.

Each microorganism is characterized by a volume \( \vartheta \) and a density contrast \( \Delta \rho_f \) relative to the fluid density \( \rho_f \). The mean velocity of the suspension is represented by \( \boldsymbol{u}^* \). The fundamental equations governing this system are given as follows:

Continuity equation
\begin{equation}\label{1}
\boldsymbol{\nabla}^* \cdot \boldsymbol{u}^* = 0.
\end{equation}

Momentum equation
\begin{align}\label{2}
\nonumber\frac{\rho_f}{\phi}\left(\frac{\partial\boldsymbol{u}^*}{\partial t^*}+\frac{1}{\phi}\boldsymbol{u}^*\cdot\boldsymbol{\nabla}^*\boldsymbol{u}^*\right)= -\boldsymbol{\nabla}^* P_e^* + \mu_f {\nabla}^{*2} \boldsymbol{u}^* -\frac{\mu_f}{K}\boldsymbol{u}^*\\- n^* \vartheta g \Delta \rho_f \hat{\boldsymbol{z}^*} - \rho_f(1-\beta (T^* - T_0))g \hat{\boldsymbol{z}^*}.
\end{align}

Here, \( P_e^* \), \( \mu_f \), \( \beta \), and \( \hat{\boldsymbol{z}^*} \) represent the hydrostatic pressure, viscosity, thermal expansion coefficient, and unit vector in the vertical direction. Also $\phi$ and $K$ denote the porosity and permeability of the medium.

Cell conservation equation
\begin{equation}\label{3}
\frac{\partial n^*}{\partial t^*} = -\boldsymbol{\nabla}^*\cdot \boldsymbol{J}^*,
\end{equation}
where the total flux is
\begin{equation*}
\boldsymbol{J}^* = (\boldsymbol{u}^* + W_c^* \langle \boldsymbol{p} \rangle) n^* - \boldsymbol{D}_B \cdot \boldsymbol{\nabla}^* n^*.
\end{equation*}

The thermal energy equation
\begin{equation}\label{4}
\varphi\left(\frac{\partial T^*}{\partial t^*} + \boldsymbol{\nabla}^* \cdot (\boldsymbol{u}^* T^*)\right) = \alpha {\boldsymbol{\nabla}}^{*2} T^*.
\end{equation}
where the thermal conductivity and the volumetric heat capacity of
water are denoted by $\alpha$ and $\varphi$, respectively.

The cell flux formulation is based on the assumption that microorganisms are purely phototactic, neglecting viscous torque effects in horizontal swimming. The diffusion tensor is considered isotropic, \( \boldsymbol{D}_B = D_B \boldsymbol{I} \)~\cite{10ghorai2005}.

{Light Intensity Equation:}
\begin{equation}\label{5}
G(\boldsymbol{x}^*, t^*) = \mathrm{I_t} \exp\left(-\kappa^* \int_H^{\boldsymbol{r}^*} n^*(\boldsymbol{r}^*, t^*) \ud\boldsymbol{r}^*\right).
\end{equation}

Here, \( \boldsymbol{r}^* \) is the vector from the microorganism to the light source, \( \mathrm{I_t} \) is the incident light intensity, \( n^* \) represents the microorganism concentration, and \( \kappa^* \) is the extinction coefficient. This formulation is based on the Lambert-Beer law, considering only absorption by microorganisms while neglecting scattering.

%%%%%%%%%%%%%%%%%%%%%%%%%%%%%%%%%%%%%%%%%%%%%%%%%%%%%%%%%	
\subsection{SWIMMING ORIENTATION}
The mean swimming direction of cells is represented by \( \langle \boldsymbol{p} \rangle \). For many microorganism species, the swimming speed remains independent of factors such as illumination, position, time, and direction~\cite{17hill1997}. The average swimming speed is denoted by \( W_c^* \). Consequently, the mean velocity of the swimming cells is given by:
\begin{equation*}
\boldsymbol{W}_c^* = W_c^* \langle \boldsymbol{p} \rangle.
\end{equation*}

The mean swimming direction, \( \langle \boldsymbol{p} \rangle \), is expressed as:
\begin{equation}\label{6}
\langle \boldsymbol{p} \rangle = M(G) \hat{\boldsymbol{z}^*},
\end{equation}
where \( \hat{\boldsymbol{z}^*} \) is the unit vector along the vertical \( Z^* \)-axis, and \( M(G) \) is the phototaxis function, defined as:
\begin{equation*}
M(G) =
\begin{cases} 
\geq 0, & G(\boldsymbol{x}^*) \leq G_c, \\
< 0, & G(\boldsymbol{x}^*) > G_c.
\end{cases}
\end{equation*}
The mean swimming direction becomes zero when \( G = G_c \) or when \( M(G) = 0 \).

For an algal suspension exposed to uniform illumination, the light intensity is influenced by the concentration of cells, as those closer to the light source cast shadows on those farther away. Thus, the light intensity \( G(\boldsymbol{x}^*, t^*) \) at a given position \( \boldsymbol{x}^* = (x^*, y^*, z^*) \) within the suspension is determined by:

\begin{equation}\label{7}
G(\boldsymbol{x}^*, t^*) = \mathrm{I_t} \exp\left(-\kappa \int_{z^*}^{H} n^*(\boldsymbol{x}', t^*) \ud\boldsymbol{z}'\right).
\end{equation}

%%%%%%%%%%%%%%%%%%%%%%%%%%%%%%%%%%%%%%%%%%%%%%%%%%%%%%%%%%%%%%%%%%%%%%%%
\subsection{BOUNDARY CONDITIONS}
In this analysis, results will be presented for both free and rigid upper horizontal boundaries. In experimental settings, the lower boundary is typically rigid, while the upper boundary may be either free or rigid.

The boundary conditions governing the system are given by
\begin{equation}\label{8}
\boldsymbol{u}^* \cdot \hat{\boldsymbol{z}^*} = 0, \quad \text{on} \quad z^* = 0, H,
\end{equation}		
\begin{equation}\label{9}
\boldsymbol{J}^* \cdot \hat{\boldsymbol{z}^*} = 0, \quad \text{on} \quad z^* = 0, H.
\end{equation}

For a rigid boundary
\begin{equation}\label{10}
\boldsymbol{u}^* \times \hat{\boldsymbol{z}^*} = 0, \quad \text{on} \quad z^* = 0, H.
\end{equation}

For a free boundary
\begin{equation}\label{11}
\frac{\partial^2}{\partial z^{*2}} (\boldsymbol{u}^* \cdot \hat{\boldsymbol{z}^*}) = 0, \quad \text{on} \quad z^* = 0, H.
\end{equation}

Also, the thermal boundary conditions
\begin{equation}\label{12}
T^* = T_0 - \Delta T , \quad \text{on} \quad z^*= 0,
\end{equation}
\begin{equation}\label{13}
T^* = T_0 , \quad \text{on} \quad z^* = H.
\end{equation}

%%%%%%%%%%%%%%%%%%%%%%%%%%%%%%%%%%%%%%%%%%%%%%%%%%%%%%%%%%%%%%%%%%%%%%%%%%
\subsection{SCALING OF THE EQUATIONS}
The equations have been scaled following the parameters used by Ghorai \& Hill~\cite{10ghorai2005} and Kuznetsov~\cite{27kuznetsov2005thermo}. The both dimensional and nondimensional variables are shown by same symbol. The governing bioconvection equations in dimensionless form are	
\begin{equation}\label{14}
\boldsymbol{\nabla} \cdot \boldsymbol{u} = 0,
\end{equation}	
\begin{align}\label{15}
\nonumber \frac{Pr^{-1}}{\phi} \left(\frac{\partial \boldsymbol{u}}{\partial t} + \frac{1}{\phi} \boldsymbol{u} \cdot \boldsymbol{\nabla} \boldsymbol{u} \right) = -\nabla P_{e} + \nabla^{2} \boldsymbol{u} - Da^{-1} \boldsymbol{u} \\
- R_B n \hat{\boldsymbol{z}} - R \hat{\boldsymbol{z}} - R_T T \hat{\boldsymbol{z}},
\end{align}	
\begin{equation}\label{16}
\frac{\partial n}{\partial t} = -\boldsymbol{\nabla} \cdot \boldsymbol{J},
\end{equation}
where
\begin{equation*}
\boldsymbol{J} = n \boldsymbol{u} + \frac{1}{Le} n V_c \langle \boldsymbol{p} \rangle - \frac{1}{Le} \boldsymbol{\nabla} n,
\end{equation*}
and 	
\begin{equation}\label{17}
\frac{\partial T}{\partial t} + \boldsymbol{\nabla} \cdot (\boldsymbol{u} T) = \boldsymbol{\nabla}^2 T.
\end{equation}

The non-dimensional parameters considered in this study include the Prandtl number, \( Pr = \mu_f /\rho_f D_T \) and the dimensionless swimming speed is given by \( V_c = W_c^* H / D_B \). The bio-convective Rayleigh number is defined as \( R_B = \bar{n^*} \vartheta g \Delta \rho_f H^{3} / \mu_f D_T \) and the thermal Rayleigh number is given by \( R_T = g \beta \Delta T^* H^{3} / \mu_f D_T \). The basic density Rayleigh number is \( R = \rho_f g H^3 / \mu_f D_T \) and the Lewis number is defined as \( Le = D_T / D_B \). Additionally, the Darcy number, which characterizes the permeability of the porous medium, is given by \( Da = K / H^2 \).

The dimensionless form of the light intensity \( G(\boldsymbol{x}, t) \) at position \( \boldsymbol{x} = (x, y, z) \) is	
\begin{equation}\label{18}
G(\boldsymbol{x}, t) = \mathrm{I_t} \exp \left(-\tau_H \int_z^1 n(\boldsymbol{x}', t) \ud\boldsymbol{x}' \right),
\end{equation} 
where \( \tau_H = \kappa^* \bar{n^*} H \) is the non-dimensional extinction coefficient.

After scaling, the boundary conditions are given by	
\begin{equation}\label{19}
\boldsymbol{u} \cdot \hat{\boldsymbol{z}} = 0, \quad \text{on} \quad z = 0,1,
\end{equation}			
\begin{equation}\label{20}
\boldsymbol{J} \cdot \hat{\boldsymbol{z}} = 0, \quad \text{on} \quad z = 0,1.
\end{equation}

For a rigid boundary	
\begin{equation}\label{21}
\boldsymbol{u} \times \hat{\boldsymbol{z}} = 0, \quad \text{on} \quad z = 0,1,
\end{equation}	
whereas for a free boundary	
\begin{equation}\label{22}
\frac{\partial^2}{\partial z^2} (\boldsymbol{u} \cdot \hat{\boldsymbol{z}}) = 0, \quad \text{on} \quad z = 0,1.
\end{equation}

Also, the thermal boundary conditions	
\begin{equation}\label{23}
T = 0, \quad \text{on} \quad z = 0,
\end{equation}	
\begin{equation}\label{24}
T = 1, \quad \text{on} \quad z = 1.
\end{equation}

%%%%%%%%%%%%%%%%%%%%%%%%%%%%%%%%%%%%%%%%%%%%%%%%%%%%%%%%%%%%%%%%%%%%%%%%%%%	

\section{The steady solution}

The steady-state solutions of Eqs.~(\ref{1})–(\ref{4}) satisfy the following conditions

\begin{equation}\label{25}
\begin{pmatrix}
\boldsymbol{u} \\ n \\ T \\ \langle \boldsymbol{p} \rangle
\end{pmatrix}
=
\begin{pmatrix}
0 \\ n_p \\ T_p \\ \langle \boldsymbol{p}_p \rangle
\end{pmatrix}.
\end{equation}

The steady-state concentration \( n_p(z) \) satisfies the equation

\begin{equation}\label{26}
\frac{\ud n_p}{\ud z} - V_c M_p n_p = 0,
\end{equation}

where \( M_p = M(G_p) \) represents the phototaxis function evaluated at steady-state light intensity \( G_p \), given by

\begin{equation}\label{27}
G_p(\boldsymbol{z}) = \mathrm{I_t} \exp \left(-\tau_H  \int_z^1 n_p(\boldsymbol{z}') \ud\boldsymbol{z}' \right).
\end{equation}    

Introducing the transformation \( \psi = \int_1^z n_p(\boldsymbol{z}') \ud\boldsymbol{z}' \), Eq.~(\ref{27}) is rewritten as

\begin{equation}\label{28}
\frac{\ud^2\psi}{\ud z^2} - V_c M_p \frac{\ud\psi}{\ud z} = 0.
\end{equation}

The boundary conditions for \( \psi \) are given by

\begin{equation}\label{29}
\psi + 1 = 0, \quad \text{at} \quad z = 0,
\end{equation}
\begin{equation}\label{30}
\psi = 0, \quad \text{at} \quad z = 1.
\end{equation}

The temperature profile \( T_p(z) \) is governed by

\begin{equation}\label{31}
\frac{\ud^2 T_p}{\ud \,z^2} = 0.
\end{equation}

Applying the boundary conditions (\ref{23}) and (\ref{24}), the solution for \( T_p(z) \) is

\begin{equation}\label{32}
T_p(z) = z.
\end{equation}

Next, Eq.~(\ref{28}) along with boundary conditions (\ref{29}) and (\ref{30}), define a boundary value problem which is solved by using a shooting method.

In this study, the phototaxis function is modeled as

\begin{equation}\label{33}
M(G) = 0.8 \sin \left[\frac{3\pi}{2} \Xi(G) \right] - 0.1 \sin \left[\frac{\pi}{2} \Xi(G) \right],
\end{equation}

where \( \Xi(G) = G \exp[\beta(G - 1)] \). The parameter \( \beta \) is related to \( G_c \), and its specific form depends on the type of microorganisms considered.
%%%%%%%%%%%%%%%%%%%%%%%%%%%%%%%%%%%%%%%%%%%%%%%%%%%%%%%%%%%%%%%%%%%%%%%	   

\section{The perturbed system}
To analyze the stability of the system, a small perturbation of amplitude \( \epsilon \) (where \( 0 < \epsilon \ll 1 \)) is introduced to the steady-state solution, leading to	
\begin{equation}\label{34}
\begin{pmatrix}
\boldsymbol{u} \\ n \\ T \\ \langle \boldsymbol{p} \rangle
\end{pmatrix}
=
\begin{pmatrix}
0 + \epsilon \boldsymbol{u}_1 + O(\epsilon^2) \\ n_p + \epsilon n_1 + O(\epsilon^2) \\ T_p + \epsilon T_1 + O(\epsilon^2) \\ \langle \boldsymbol{p}_p \rangle + \epsilon \langle \boldsymbol{p}_1 \rangle + O(\epsilon^2)
\end{pmatrix},
\end{equation}	
where \( \boldsymbol{u}_1 = (u_1, v_1, w_1) \). 

Substituting these perturbed variables into Eqs.~(\ref{1})–(\ref{4}) and retaining only terms of order \( O(\epsilon) \) results in the following linearized system	
\begin{equation}\label{35}
\boldsymbol{\nabla} \cdot \boldsymbol{u}_1 = 0,
\end{equation}	
\begin{equation}\label{36}
\frac{Pr^{-1}}{\phi} \left( \frac{\partial \boldsymbol{u}_1}{\partial t} \right) = -\boldsymbol{\nabla} P_{e} + \nabla^{2} \boldsymbol{u}_1 - Da^{-1} \boldsymbol{u}_1 - R_B n_1 \hat{\boldsymbol{z}} - R_T T_1 \hat{\boldsymbol{z}},
\end{equation}	
\begin{equation}\label{37}
\frac{\partial n_1}{\partial t} + \frac{1}{Le} V_c \boldsymbol{\nabla} \cdot \left( \langle \boldsymbol{p}_p \rangle n_1 + \langle \boldsymbol{p_1} \rangle n_p \right) + w_1 \frac{\ud n_p}{\ud\,z} = \frac{1}{Le} \boldsymbol{\nabla}^2 n_1,
\end{equation}	
\begin{equation}\label{38}
\frac{\partial T_1}{\partial t} + w_1 \frac{\ud T_p}{\ud\,z} = \boldsymbol{\nabla}^2 T_1.
\end{equation}

If \( G = G_p + \epsilon G_1 \), the perturbed light intensity \( G_1 \), after expanding \( \exp(-\tau_H  \int_z^1 \epsilon n_1 dz') \) and retaining terms of order \( O(\epsilon) \), is given by	
\begin{equation}\label{39}
G_1(\boldsymbol{z}) =\tau_H G_p \left( - \int_z^1  n_1(\boldsymbol{z}') \ud \boldsymbol{z}' \right).
\end{equation}

Expanding \( M(G_p + \epsilon G_1) \hat{\boldsymbol{z}} - \langle \boldsymbol{p}_p \rangle \) to first order in \( \epsilon \) gives	
\begin{equation}\label{40}
\langle \boldsymbol{p}_1 \rangle = G_1 \frac{dM_p}{dG} \hat{\boldsymbol{z}}.
\end{equation}

Substituting \( \boldsymbol{p}_1 \) into Eq.~(\ref{37}) and simplifying	
\begin{align}\label{41}
\nonumber \frac{\partial n_1}{\partial t} + w_1 \frac{\ud n_p}{\ud\,z} + \frac{1}{Le} \aleph_1(z) \int_1^z n_1 dz' + \frac{1}{Le} \aleph_2(z) n_1 \\
+ \frac{1}{Le} \aleph_3(z) \frac{\partial n_1}{\partial z} = \frac{1}{Le} \nabla^2 n_1,
\end{align}	
where
\begin{equation*}
\aleph_1(z) = \tau_H  V_c \frac{\ud}{\ud z} \left( n_p G_p \frac{\ud M_p}{\ud G} \right),
\end{equation*}
\begin{equation*}
\aleph_2(z) = 2\tau_H  V_c n_p G_p \frac{\ud M_p}{\ud G},
\end{equation*}
\begin{equation*}
\aleph_3(z) = V_c M_p.
\end{equation*}

To eliminate \( P_e \) and the horizontal components of \( \boldsymbol{u}_1 \), we take the double curl and extract the vertical component, reducing Eqs.~(\ref{35})-- (\ref{38}) to three equations for \( w_1 \), \( n_1 \), and \( T_1 \). These are then expressed in normal modes as	
\begin{equation}\label{42}
\begin{pmatrix}
w_1 \\ n_1 \\ T_1
\end{pmatrix}
=
\begin{pmatrix}
W(z) \\ \Theta(z) \\ T(z)
\end{pmatrix}
\exp{[\gamma t + i(lx + my)]}.
\end{equation}

Thus, the governing equations simplify to	
\begin{align}\label{43}
\nonumber \left( \gamma \frac{Pr^{-1}}{\phi} + Da^{-1} + k^2 - \frac{d^2}{dz^2} \right) \left( \frac{d^2}{dz^2} - k^2 \right) W \\
= R_B k^2 \Theta(z) + R_T k^2 T(z),
\end{align}	
\begin{align}\label{44}
\nonumber D^2 \Theta - \aleph_3(z) D \Theta - (Le\gamma + k^2 \aleph_2(z)) \Theta - \aleph_1(z) \int_1^z \Theta dz \\
- \aleph_0(z) = Le D n_p W,  
\end{align}	
\begin{equation}\label{45}
\left( \frac{d^2}{dz^2}-\gamma - k^2   \right) T(z) = W(z).
\end{equation}  

The boundary conditions are	
\begin{align}\label{46}
\nonumber W = \frac{dW}{dz} = \frac{d\Theta}{dz} - \aleph_3(z) \Theta + n_p V_c \tau_H \\
\times G_p \left( \int_z^1 \Theta(z') dz' \right) \frac{dM_p}{dG} = 0, \quad \text{on} \quad z = 0,
\end{align}	
\begin{equation}\label{47}
W = \frac{d^2 W}{dz^2} = \frac{d\Theta}{dz} - \aleph_3(z) \Theta = 0, \quad \text{on} \quad z = 1.
\end{equation}

For a rigid upper boundary, the condition in Eq.~(\ref{47}) is replaced by	
\begin{equation}\label{48}
W = \frac{dW}{dz} = \frac{d\Theta}{dz} - \aleph_3(z) \Theta = 0, \quad \text{on} \quad z = 1.
\end{equation}

The thermal boundary conditions are	
\begin{equation}\label{49}
T(z) = 0, \quad \text{on} \quad z = 0,1.
\end{equation}

The system Eqs.~(\ref{43})--(\ref{45}) form an eigenvalue problem for \( \gamma \), determining the stability of the perturbations.

Introducing a new variable:	
\begin{equation*}
\Phi(z) = \int_1^z \Theta(z') dz'.
\end{equation*}

The perturbed equations, rewritten using \( D = d/dz \), take the form	
\begin{align}\label{50}
\left( D^2 - k^2 \right) \left( \gamma \frac{Pr^{-1}}{\phi}+Da^{-1} + k^2 - D^2 \right) W = R_B k^2 D \Phi + R_T k^2 T.
\end{align}	

\begin{align}\label{51}
\nonumber D^3 \Phi - \aleph_3(z) D^2 \Phi - (\gamma Le + k^2 + \aleph_2(z)) D \Phi - \aleph_1(z) \Phi \\
- \aleph_0(z) = Le D n_p W,
\end{align}	
\begin{equation}\label{52}
\left( D^2-k^2 - \gamma \right) T =  W.
\end{equation}

The boundary conditions are	
\begin{align}\label{53}
W = DW = \aleph_2(z) \Phi + 2 \aleph_3(z) D \Phi - 2 D^2 \Phi = 0, \quad \text{on} \quad z = 0.
\end{align}		
\begin{align}\label{54}
W = D^2 W = \aleph_2(z) \Phi + 2 \aleph_3(z) D \Phi - 2 D^2 \Phi = 0, \quad \text{on} \quad z = 1.
\end{align}	

For a rigid upper surface, the boundary condition in Eq.~(\ref{54}) is replaced by	
\begin{align}\label{55}
W = DW = \aleph_2(z) \Phi + 2 \aleph_3(z) D \Phi - 2 D^2 \Phi = 0, \quad \text{on} \quad z = 1.
\end{align}	
\begin{equation}\label{56}
T = 0, \quad \text{on} \quad z = 0,1.
\end{equation}	

Additionally, an extra boundary condition	
\begin{equation}\label{57}
\Phi = 0, \quad \text{on} \quad z = 1.
\end{equation}

%%%%%%%%%%%%%%%%%%%%%%%%%%%%%%%%%%%%%%%%%%%%%%%%%%%%%%%%%%%%%%%%%%%%%%%%%%%

\section{SOLUTION technique}
The stability analysis of the system described by Eqs.~(\ref{50})–(\ref{52}) is carried out using a fourth-order finite difference scheme combined with Newton-Raphson-Kantorovich (NRK) iterations~\cite{18cash1980}. This numerical technique allows for the evaluation of the growth rate, Re$(\gamma)$, and the construction of neutral stability curves in the $(k, R)$-plane under different parameter conditions, where \( R \) can take values \( R_B \) or \( R_T \).

 The neutral stability curve, expressed as \( R^{(n)}(k) \), contains an infinite sequence of branches indexed by \( n = 1, 2, 3, \dots \), each corresponding to a distinct solution of the linear stability problem. The most unstable mode is associated with the lowest \( R \) value, defining the critical solution as \( (k_c, R^c) \), where the critical Rayleigh number \( R^c \) corresponds to either \( R_B^c \) or \( R_T^c \).

%%%%%%%%%%%%%%%%%%%%%%%%%%%%%%%%%%%%%%%%%%%%%%%%%%%%%%%%%%%%%%%%%%%%%	
\section{NUMERICAL RESULTS}

This study systematically investigates the influence of thermal effects on algal suspensions by exploring the parameter space. A controlled approach is adopted, where specific parameters remain fixed while others are varied to assess their impact on bioconvection onset. To maintain consistency, the governing parameters are set as \( Pr = 5 \) and \( \mathrm{I_t} = 0.8 \), while different values of other parameters are considered to analyse their effects.  

To account for variations in light absorption properties, two cases are examined: \( \tau_H = 0.5 \) and \( \tau_H = 1.0 \). Furthermore, the critical light intensity, \( G_c \), is chosen such that the steady-state condition occurs at different depths within the suspension, specifically at the top, three-quarters, or mid-height. The primary objective of this analysis is to determine the critical bioconvective Rayleigh number, \( R_B^c \), as a function of the thermal Rayleigh number, \( R_T \), providing insights into the stability conditions under varying thermal influences.

%%%%%%%%%%%%%%%%%%%%%%%%%%%%%%%%%%%%%%%%%%%%%%%%%%%%%%%%%%%%%%%%%%%%%%%%%%%%
\subsection{$V_c=10$}
The results of this section is divided into two categories based on the top surface conditions.

\subsubsection{WHEN TOP SURFACE IS STRESS FREE}
\begin{figure}[!htbp]
	\centering
	\includegraphics[scale=0.5]{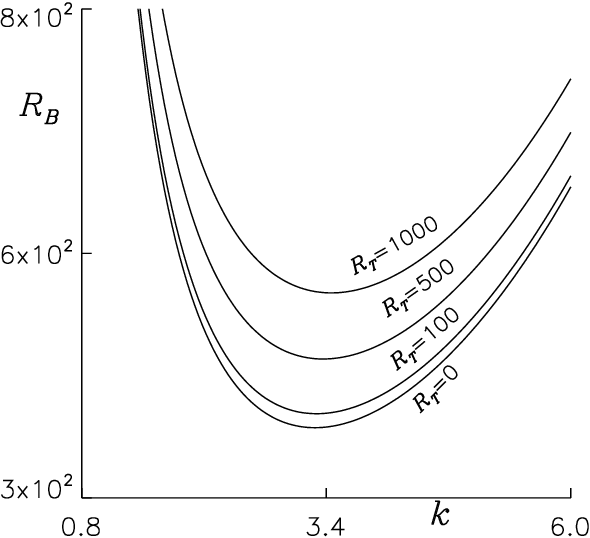}
	\caption{\footnotesize{Neutral curves for different values of the thermal Rayleigh number, where the other parameters $V_c=10,\tau_H=0.5,G_c=0.63,\phi=0.76,Da=0.1$, and $Le=0.4$ are fixed. Here, the top wall is assumed to be stress free.}}
	\label{fig2}
\end{figure}

Figure~\ref{fig2} shows the neutral curves for different values of the thermal Rayleigh numbers. Here the other parameters $V_c=10,\tau_H=0.5,G_c=0.63,\phi=0.76,Da=0.1$, and $Le=0.4$ are fixed. It is clear that the as the thermal Rayleigh number increases, the critical bioconvective Rayleigh number increases. Thus the suspension becomes more stable as the thermal Rayleigh number increases.

\begin{figure}[!htbp]
	\centering
	\includegraphics[scale=0.5]{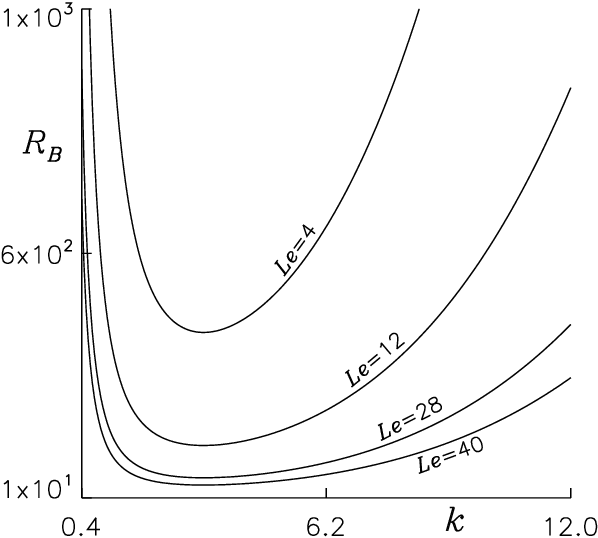}
	\caption{\footnotesize{Neutral curves for different values of the Lewis number, where the other parameters $V_c=10,\tau_H=0.5,G_c=0.63,\phi=0.76,Da=0.1$, and $R_T=50$ are fixed. Here, the top wall is assumed to be stress free.}}
	\label{fig3}
\end{figure}

Also, Fig.~\ref{fig3} illustrates the neutral curves for different values of the Lewis numbers. Here the other parameters $V_c=10,\tau_H=0.5,G_c=0.63,\phi=0.76,Da=0.1$, and $R_T=50$ are fixed. The critical bioconvective Rayleigh number decreases as the thermal Rayleigh number increases.

%%%%%%%%%%%%%%%%%%%%%%%%%%%%%%%%%%%%%%%%%%%%%%%%%%%%%%%%%%%%%%%%%%%%%%
\subsubsection{WHEN TOP SURFACE IS RIGID}
\begin{figure}[!htbp]
	\centering
	\includegraphics[scale=0.5]{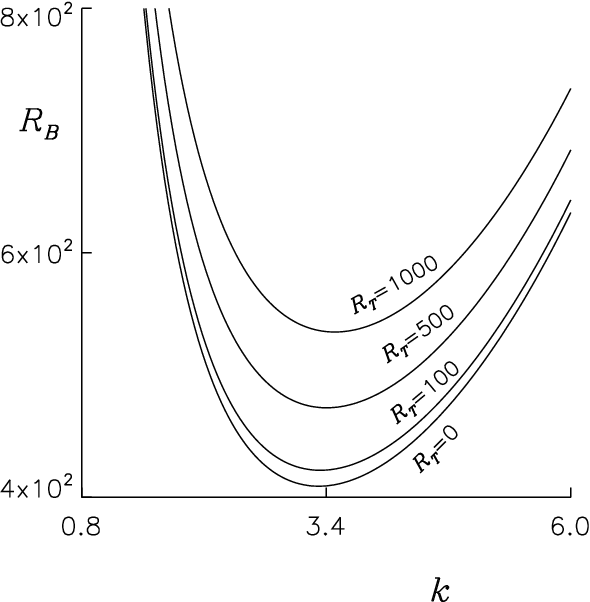}
	\caption{\footnotesize{Neutral curves for different values of the thermal Rayleigh number, where the other parameters $V_c=10,\tau_H=0.5,G_c=0.63,\phi=0.76,Da=0.1$, and $Le=0.4$ are fixed. Here, the top wall is assumed to be rigid.}}
	\label{fig4}
\end{figure}

Figure~\ref{fig4} shows the neutral curves for different values of the thermal Rayleigh numbers. Here the other parameters $V_c=10,\tau_H=0.5,G_c=0.63,\phi=0.76,Da=0.1$, and $Le=0.4$ are fixed. It is clear that the as the thermal Rayleigh number increases, the critical bioconvective Rayleigh number increases. Thus the suspension becomes more stable as the thermal Rayleigh number increases.

Also, Fig.~\ref{fig3} illustrates the neutral curves for different values of the Lewis numbers. Here the other parameters $V_c=10,\tau_H=0.5,G_c=0.63,\phi=0.76,Da=0.1$, and $R_T=50$ are fixed. The critical bioconvective Rayleigh number decreases as the thermal Rayleigh number increases.

\begin{figure}[!htbp]
	\centering
	\includegraphics[scale=0.5]{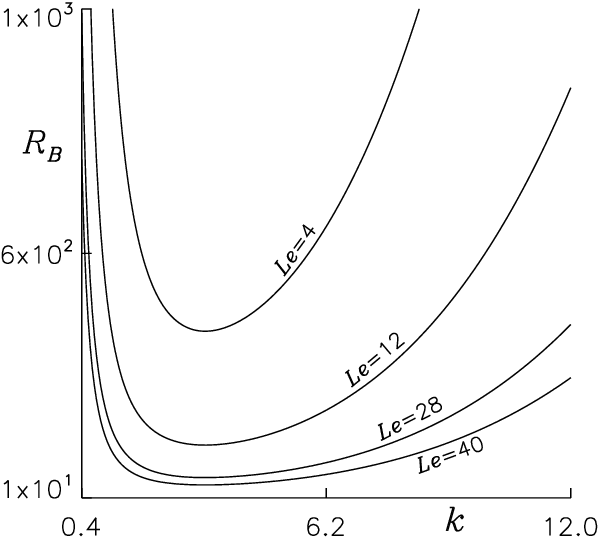}
	\caption{\footnotesize{Neutral curves for different values of the Lewis number, where the other parameters $V_c=10,\tau_H=0.5,G_c=0.63,\phi=0.76,Da=0.1$, and $R_T=50$ are fixed. Here, the top wall is assumed to be rigid.}}
	\label{fig5}
\end{figure}

%%%%%%%%%%%%%%%%%%%%%%%%%%%%%%%%%%%%%%%%%%%%%%%%%%%%%%%%%%%%%
\section{Conclusion}

This study presents a novel thermal phototactic bioconvection model that integrates thermal effects into the dynamics of bioconvective systems in a porous medium. The suspension is exposed to collimated light from above and also subjected to heating from the same side. The primary objective is to examine how thermal variations impact the onset and behavior of phototactic bioconvection in a porous environment. The numerical results, obtained through a linear stability analysis under specific control parameters, are summarized below.

The accumulation of algal cells within the suspension forms a sublayer, whose position is governed by the critical light intensity. When the critical intensity matches the incident light intensity, cells cluster at the top; as the critical intensity decreases, this aggregation shifts downward. Moreover, a reduction in light intensity corresponds to a decrease in the maximum concentration at the sublayer. 

The linear stability analysis reveals both steady and oscillatory solutions, with oscillations becoming particularly pronounced when the microorganism layer is positioned around three-quarters of the suspension height, dictated by the critical light intensity. Notably, these oscillations transition to a stationary state as the critical thermal Rayleigh number increases. As the thermal Rayleigh number increases, the critical bioconvective Rayleigh number increase. On the other hand, the critical bioconvective Rayleigh number decreases as the Lewis number increases. Thus the suspension become stable (unstable) as the thermal Rayleigh number (Lewis number) increases.  

%%%%%%%%%%%%%%%%%%%%%%%%%%%%-OM NAMAH SHIVAY-%%%%%%%%%%%	
%\section*{AUTHOR DECLARATIONS}
\section*{CONFLICT OF INTEREST}
The authors state that they have no conflicts of interest.
\section*{DATA AVAILABILITY}
The findings of this study are supported by the data contained within this article.
\nocite{*}
\section*{REFERENCES}
\bibliography{COLLIMATED_HEATED_FROM_TOP}

\end{document}